%Paper: q-alg/9505014
%From: fronsdal <fronsdal@physics.ucla.edu>
%Date: Tue, 16 May 95 15:05:48 0800

%\\
%Title: Universal T-matrix for Twisted Quantum gl(N).
% Author: Christian Fronsdal, UCLA.
%Comments: 30 pages. Plain TeX.
%\\
%Abstract. The Universal T-matrix is the capstone of the structure that
%consists of a quantum group and its dual, and the central object from which
%spring the T-matrices (monodromies) of all the associated integrable models.  A
%closed expression is obtained for the case of multiparameter (twisted) quantum
%$gl(N)$.  The factorized nature of standard quantum groups, that allows the
%explicit expression for $U\hskip-1mm T $ to be obtained with relative ease,
%extends to some nonstandard quantum groups, such as those based on $A_n^{(2)}$,
%and perhaps to all.  The paper is mostly concerned with parameters in general
%position, but the extension to roots of unity is also explored, in the case of
%$g\ell(N)$.  The structure of the dual is now radically different, and an
%interesting generalization of the $q$-exponential appears in the formulas for
%the Universal T- and R-matrices. The projection to quantum $sl(N)$ is simple
%and
%direct; this allows, in particular, to apply recent results concerning
%deformations of twisted $gl(N)$ to the semisimple quotient.
%\\

%Beginning Formats

%\magnification=1200
\magnification=\magstep1
\vsize=21.5truecm
\voffset=1.00truein
\settabs 18 \columns
%\hoffset=3.75truecm
\hoffset=1.00truein
\hsize=15.2truecm
%\nopagenumbers
\baselineskip=17 pt

\font\steptwo=cmr10 scaled\magstep2

\def\b{\bigskip}
\def\bb{\bigskip\bigskip}
\def\bbb{\bigskip\bigskip\bigskip}

\def\no{\noindent}

\def\ce{\centerline}
\def\ve{\vfill\eject}

\def\e e{$e^+ e^-$ }

\def\A{$A(<\hskip-1mm q \hskip-1mm > ,a)\,\,$}
%End of Beginning Formats
%Beginning of Letter Heading

\rightline {UCLA/93/TEP/3}
\rightline {February 1993}

\bbb
\ce {\steptwo {Universal T-matrix for Twisted Quantum $gl(N)$}}
\bb
\ce {\bf C. Fr{\o}nsdal}
%\vskip1.50truein
\b
\ce {\it Department of Physics}
\ce {\it University of California, Los Angeles, CA 90024-1547}

\bbb
\ce {\bf Abstract}
\b The Universal T-matrix is the capstone of the structure that consists of a
quantum group and its dual, and the central object from which spring the
T-matrices (monodromies) of all the associated integrable models.  A closed
expression is obtained for the case of multiparameter (twisted) quantum
$g\ell(N)$.  The factorized nature of standard quantum groups, that allows the
explicit expression for $U\hskip-1mm T $ to be obtained with relative ease,
extends to some nonstandard quantum groups, such as those based on $A_n^{(2)}$,
and perhaps to all.  The paper is mostly concerned with parameters in general
position, but the extension to roots of unity is also explored, in the case of
$g\ell(N)$.  The structure of the dual is now radically different, and an
interesting generalization of the $q$-exponential appears in the formulas for
the Universal T- and R-matrices. The projection to quantum $sl(N)$ is simple
and
direct; this allows, in particular, to apply recent results concerning
deformations of twisted $gl(N)$ to the semisimple quotient.

\ve

\ce {\bf 1. Introduction}
\b
\no  1.1. PICTURES

Quantum groups may be studied in two different ``pictures", much as quantum
mechanics has a Schroedinger picture and a Heisenberg picture.  The Drinfeld
picture is a quantum group from the point of view of integrable field theories,
where the first examples were discovered.  The development is due to Drinfeld
[1], Jimbo [2] and many others.  The Woronowicz-Manin picture appeared already
in early work of Baxter [3] and Faddeev {\it et al.} [4], but the systematic
study is due to Woronowicz [5] and Manin [6].  The term psuedo-group, proposed
by Woronowicz in [5], deserves to be retained, in order that the terminology
distinguish between this object and the quantum group in the sense of Drinfeld.

\vskip.5cm

\no 1.2. DUALITY

This paper is a contribution to the study of the duality between the two
pictures.  Much work has been done already [7], and a firm foundation that
incorporatess reflexivity and deformation theory has been discovered very
recently [8], but the object that we regard as the capstone of the structure
seems not to have been calculated except in some simple cases.  It appears as
the universal bi-character in Woronowicz and as a ``canonical element" or
``dual
form" in other places.  Here it will be called the Universal T-matrix, in
recognition of the fact that the transition matrices of integrable models
appear
upon specialization, in passing from structure to representations.  This is the
main reason why it is useful for physics.

Let $A$ be a vector space with a (countable) basis $\{|\alpha>\}$.  Let
$\{<\beta|\}$ be the functions on $A$ defined by

$$<\beta|\alpha>=\delta_{\alpha\beta}.\eqno(1.1)$$

\no To avoid or at least postpone the thorny question of an appropriate and
precise definition of the dual space $A^*$, let us just say that these
functions, and (at least) all finite linear combinations of them, are in $A^*$.
The Universal T-matrix is the unit operator in $End\,\, A$, namely

$$U\hskip-1mm T = \sum_\alpha |\alpha><\alpha|,\eqno(1.2)$$

\no the ``resolution of the identity."  It has been called the  ``canonical
element of ${A \otimes A^*}$" and this expression can perhaps be justified once
a
convenient and rigorous definition of $A^*$ has been agreed upon.

The notation was chosen in order to exhibit the connection between the object
of
main interest and familiar operations in quantum mechanics in Eq. (1.2).  We
now
switch to a more convenient notation.

Let $A$ be an associative algebra with basis $\{L_\alpha\}$, and $\{F_\alpha\}$
the functions defined by

$$F_\alpha(L_\beta)=\delta_{\alpha\beta}.\eqno(1.3)$$

\no We shall suppose that the product $L_\alpha L_\beta$ is a finite linear
combination of the basis elements; then a structure of coassociative coproduct
is naturally induced on the linear span $A^*$ of the functions $F_\alpha$; it
is
defined by

$$F_\alpha \mapsto \Delta F_\alpha,\quad \Delta F_\alpha(L_\beta,L_\gamma)
   =F_\alpha(L_\beta L_\gamma).\eqno(1.4)$$

\no The expression for $U\hskip-1mm T$ is

$$U\hskip-1mm T = T_{_{L,F}} = \sum_\alpha L_\alpha F_\alpha,\eqno(1.5)$$

\no and the relation (1.4) can be expressed as

$$T_{_{L,F}} T_{_{L,F^\prime}}=\sum_{\alpha\beta}
  L_\alpha L_\beta F_\alpha F_\beta^\prime = \sum_\gamma
  L_\gamma \Delta F_\gamma.\eqno(1.6)$$

\no Here $F,F^\prime$ stands for two copies of $A$, $F_\alpha F_\beta^\prime$
has the same meaning as $F_\alpha \otimes F_\beta$, and the result is that

$$T_{_{L,F}} T_{_{L,F^\prime}}=T_{_{L,\Delta F}}.\eqno(1.7)$$

\no This is the first of the two structural relations that characterize the
universal T-matrix.

Keep the notations as above, and suppose that $A$ has, in addition, a coproduct
that turns it into a bialgebra.  Then $A^*$ also becomes a bialgebra, with
algebraic structure uniquely defined by

$$F_\alpha F_\beta(L_\gamma)=F_\alpha \otimes F_\beta(\Delta L_\gamma).
  \eqno(1.8)$$

\no This too can be expressed in terms of $U\hskip-1mm T$:

$$T_{_{L,F}} T_{_{L^\prime,F}}=T_{_{\Delta L,F}}.\eqno(1.9)$$

\no This is the second structural relation satisfied by the Universal T-matrix;
note that (1.7) and (1.9) are equivalent to (1.4) and (1.8).

Finally, suppose that $A$, as an algebra, is finitely generated by elements
$\{\ell_i\}\,\,i=1,\cdots,n$, and that each basis element $L_\alpha$ is an
ordered monomial.  We would like to answer the following.
\vskip.5cm
\no {\bf Question.}  Under what condition is $A^*$ finitely generated (as an
algebra), and under what additional conditions is each element of the dual
basis
$\{F_\alpha\}$ a polynomial in the generators?
\vskip.5cm
\no Though we do not know the answer in general, we can always use the
relations
(1.7) and (1.9) to determine the structure of $A^*$ and thus answer the
question
on a case by case basis.  Furthermore, when $A^*$ turns out to be finitely
generated we expect to obtain a useful expression for
$U\hskip-1mm T$ in terms of the two sets of generators.

\vskip.5cm
\no 1.3. EXAMPLES

Let $A$ be the unital algebra freely generated by a single element $x$, with

$$\Delta x = x \otimes 1+1 \otimes x,$$

\no and take the basis $1,x,x^2,\cdots$.  Then

$$U\hskip-1mm T = \sum F_nx^n,$$

\no and the problem is to determine the dual structure.  Eq. (1.9) gives

$$\sum F_mF_n x^m \otimes x^n=\sum F_k(x \otimes 1+1\otimes x)^k$$

\no or

$$F_kF_1=(k+1)F_{k+1}\eqno(1.10)$$

\no with the unique solution ($p:=F_1$)

$$F_k=p^k/k!$$

\no Hence $A^*$ is the unital algebra freely generated by $p$, and

$$U\hskip-1mm T=e^{xp}.\eqno(1.11)$$

\no Finally, (1.7) provides an easy evaluation of the coproduct:

$$\eqalign{&e^{xp}e^{xp^\prime}=e^{x\Delta p} \cr
   &\Rightarrow \Delta p=p+p^\prime = p\otimes 1+1 \otimes p.\cr}$$

Next, let $G$ be a Lie group and ${\cal{G}}$ its Lie algebra with basis
$\{\ell_i\}\,\,i=1,\cdots,n$.  Let A be the universal enveloping algebra of
$ {\cal {G}} $, with a basis $\{L_\alpha\}$ of ordered monomials, and $\Delta$
the unique compatible coproduct generated by

$$\Delta \ell_i=\ell_i \otimes 1+1 \otimes \ell_i.\eqno(1.12)$$

\no [Compatible, that is, with the structure; $\ell_i\mapsto \Delta\ell_i$
generates a homomorphism.]  Let $\{F_\alpha\}$ be the terms of functions on $G$
at the identity defined by

$$L_\alpha F_\beta\bigm|_{Id}=\delta_{\alpha\beta}.$$

\no Then a calculation that follows precisely the pattern of the first example
leads to a unique structure of bialgebra on the space $A^*$ spanned by
$\{F_\alpha\}$.  This structure depends on the choice of the basis
$\{L_\alpha\}$.

(i) If $L_\alpha$ is the symmetric product of $\{\ell_i\}\,\,i\in \alpha$, then

$$F_\alpha = (1/|\alpha|!) \prod_{i\in\alpha} p_i.$$

\no The set $\alpha$ includes repetitions, $i\in\{1,\cdots,n\}$, and
$|\alpha |$) is the cardinality of $\alpha$.  The structure of $A^*$ is
Abelian,

$$U\hskip-1mm T = e^{p\cdot\ell},\quad p\cdot\ell:=\sum^n_{i=1} p_i\ell_i
\eqno(1.13)$$

\no and the coproduct $p_i \mapsto \Delta p_i$ is given by the
Campbell-Haussdorff formula.

(ii) If $L_\alpha$ is an ordered polynomial, defined in terms of an ordering of
the $\ell_i$, for example $\ell_1<\ell_2<\cdots < \ell_n$, then

$$F_\alpha = \mathrel{\mathop{\rm{II}}^n_{i=1}}(p_i)^{\alpha_i}/\alpha_i!,$$

\no where $\alpha_i$ is the incidence of $i$ in $\alpha$.  The structure is
Abelian,

$$U\hskip-1mm T = \prod^n_{i=1}e^{p_i\ell_i} \eqno(1.14)$$

\no with the same ordering of factors, and the coproduct is expressed by
another
version of the Campbell-Hausdorff formula.
\vskip.5cm
\no 1.4. QUANTUM GROUPS

Our program is to obtain a formula for the Universal T-matrix of quantum groups
analogous to (1.11) and (1.14).  The problem of a quantum version of (1.13)
seems to be more difficult [9].

For the quantum group $U_{q,q^\prime}(g\ell_2)$ the formula

$$U\hskip-1mm
T=e_a^{p_-x_-}e^{p_1\rho_1}e^{p_2\rho_2}e_{1/a}^{p_+x_+},\eqno(1.15)$$

$$e_a^z := \sum z^n/[n!]_a,\,\,a=1/qq^\prime,$$

\no was given in [10].  The $p_i$ generate the two-parameter (twisted) version
of quantum $g\ell(2)$ in the sense of Drinfeld, with

$$[p_-,p_+]= (q - 1/q')(q^{p_1}q'^{p_2} - q'^{-p_1}q^{-p_2}) $$

\no and the rest; the generators $x_\pm,\rho_i$ of the dual generates a
solvable
Lie algebra with

$$[\rho_i,x_\pm] \propto x_\pm ,$$

$$[\rho_1,\rho_2]=0=[x_+,x_-].$$

\no This is a special case of twisted, quantum $g\ell(N)$, treated in detail in
this paper, with analogous results for all $N$.

The key to a simple generalization of (1.15) is to choose a preferred ordering
of the generators.  In this paper we take the Woronowicz-Manin pseudogroup
($A$) as our starting point, to end up (via the construction of
$U\hskip-1mm T$) with the Drinfeld quantum group ($A^*$).  The strategy is
exactly the same as that followed in the above examples.  The important
question
of ordering is related to a unique factorization of $A$, carried out in two
steps.  First

$$A=A_- \otimes_{A_0} A_+,\eqno(1.16)$$

\no where $A_\pm$ and $A_0$ are sub-bialgebras.  Then

$$A_- \ni \tilde X = \prod^k_{i=1} \tilde X(i),\quad  A_+\ni \tilde Y=
\prod^k_{i=1} \tilde Y(i) \eqno(1.17)$$

\no where $\tilde X(i)$ and $\tilde Y(i)$ are nilpotent matrices arranged in a
particular order, the choice of which is quite crucial.

The factorization, for the case of twisted, quantum $g\ell(N)$  is carried out
with all details in Section 2, the Universal T-matrix is found in Section 3 and
the structure of the dual in Section 4.  In this case $A_0$ is Abelian and the
structure expressed by (1.16) is essentially the famous quantum double.

In Section 5 we show how the Universal R-matrix is obtained by a simple
projection from $UT$, and in Section 6 we discuss the generalizations of (1.15)
and (1.16) to
\item{(i)} The standard quantization of any simple Lie algebra.
\item{(ii)} Nonstandard quantum groups $(A_n^{(2)},D_n^{(2)})$.
\item{(iii)} Roots of unity.

Finally, in Section 7, we study the relation between twisted, quantum
$g\ell(N)$ and $s\ell(N)$, and the question of rigidity of these quantum groups
to further, essential deformations.

This paper does not include physical applications, but it may be worth
repeating
that the transfer matrices (monodromies) of solvable models (without spectral
parameters) are obtained from UT by specialization to a representation of $A$.
The choice of generators of $A$ that is introduced via the factorization (1.16)
and (1.17), gives a presentation of $A$ as a deformed enveloping algebra, which
greatly facilitates the construction of representations.  This brings out the
interesting problem of giving a sense to the Universal T-matrix for an affine
quantum group (with spectral parameter).  Applications without spectral
parameter include ice-type models, knot theory and conformal field theory in
two
dimensions, but the inclusion of a spectral parameter seems to be required for
more interesting physical applications.  All this work on finite dimensional
quantum groups must be regarded as preparation for an assault on the real
problem.
\vskip.5cm
\no 1.5. TWISTED QUANTUM $g\ell(n)$

This quantum group was found independently by Reshetikhin [11], Schirrmacher
[12]
and Sudbery [13].  In [11] it is offered as an example of ``gauge
transformations" on quantum groups based on simple Lie algebras.  Since
$g\ell(N)$ is not simple it is of some interest to study the relationship
between twisted $g\ell(N)$ and twisted $s\ell(N)$.

The R-matrix of twisted, quantum $g\ell(N)$ is an element $R\in End(V\otimes
V)$, $V=N$-dimensional vector space, that in a particular basis is given by

$$R=\sum_i M_i^i\otimes M_i^i+\sum_{i<j}
  (q^{ji}M_j^j\otimes M_i^i+aq^{ij}M_i^i\otimes M_j^j+(1-a)
  M_j^i\otimes M_i^j).\eqno(1.18)$$

\no The parameters are $q^{ij}$ with $q^{ij}q^{ji}=1$ and $a$.  The $q$'s are
subject to change under ``gauge transformation"; the parameter $a$ is much more
fundamental.  The associated deformed permutation matrix $P$, often denoted
$\check R$, is defined by

$$P_{ij}^{k\ell}=R_{ji}^{k\ell};\eqno(1.19)$$

\no it has two eigenvalues and satisfies the Hecke condition

$$(P-1)(P+a)=0.\eqno(1.20)$$

\no Reference to ``roots of unity" means that $a^K=1$ for some (minimal)
positive integer $K$.  The observation that the Hecke parameter $a$, and not
the
$q$'s, is the relevant parameter, raises the question of identifying the
important parameters in the case of quantum groups that do not satisfy the
Hecke
condition.

The matrix defined by (1.18) satisfies the Yang-Baxter relation

$$R_{12}R_{13}R_{23}=R_{23}R_{13}R_{12}\eqno(1.21)$$

\no and the matrix $P$ satisfies, in consequence, the braid relation

$$(braid)_{123}\equiv P_{12}P_{23}P_{12}-P_{23}P_{12}P_{23}=0.\eqno(1.22)$$

\no The next paragraph depends on (1.20) and (1.22), but not on the specific
form of $P$ that is implied by (1.18) and (1.19).

A quantum plane, in the sense of Manin, is an algebra generated by
$\{x^i\}\,\,i=1,\cdots,N$, with relations

$$xx(P-1)=0\quad (x^ix^jP_{ij}^{k\ell}=x^kx^\ell).\eqno(1.23)$$

\no A differential calculus $D$ on this quantum plane, as defined most
succinctly in [14], is the algebra generated by $\{x^i\}$ and
$\{\theta^i\}\,\,i=1,\cdots,N$, with additional relations

$$\theta\theta(P+a)=0,\eqno(1.24)$$

$$a\theta x=x\theta P.\eqno(1.25)$$
\vskip.5cm
\no{\bf Definition.} The pseudogroup $A(P)$ is the unital algebra generated by
the matrix elements $Z_i^j$ of an $N$-dimensional square matrix $Z$, with
relations

$$[P,Z\otimes Z]=0.\eqno(1.26)$$

The following remarks also depend only on the validity of the Hecke relation
and
the braid relation.
\vskip.5cm
\no {\bf Remark 1.}  One can look upon $A(P)$ as the algebra of quantum
automorphisms of the differential calculus $D$, the action on $D$ being
generated by $x\mapsto xZ$, $\theta\mapsto \theta Z$ ($x^i\mapsto  x^jZ_j^i$
etc.).  The idea of ``quantum automorphsims of the quantum plane" does not lead
to $A(P)$.
\vskip.5cm
\no {\bf Remark 2.}  The precise implication of the braid relation for the
differential calculus $D$ is contained in a theorem [15] that may be
paraphrased
as follows: The braid relation (1.22) $\underline{\hbox{is equivalent to}}$ the
statement ``$x$ commutes with (1.24) and $\theta$ commutes with (1.23)."
\vskip.5cm
\no {\bf Remark 3.}  An implication of (1.26) is the existence of a unique
compatible coproduct on $A(P)$, such that

$$\Delta(z_i^j)=z_i^k\otimes z_k^j.\eqno(1.27)$$

\no Henceforth we regard $A(P)$ as a bialgebra with this coproduct.
\vskip.5cm
\no {\bf Remark 4.}  The statement (1.26) can be expressed as

$$ZZ({\cal{P}}-1)=0,\eqno(1.28)$$

\no where ${\cal P}$ is the endomorphism defined by the matrix {\cal{P}} acting
on $End\,\, V\otimes V$.  This matrix does not satisfy the Hecke condition
(1.20), but instead

$$({\cal{P}}-1)({\cal{P}}+a)({\cal{P}}+a^{-1})=0.\eqno(1.29)$$

\no thus $A(P)$ is a quantum plane with ${\cal{P}}$ now taking  the place of
$P$
and (1.29) replacing (1.20).  Of course ${\cal {P}}$ satisfies the braid
relation.
\vskip.5cm

It is natural to ask what is the algebra of quantum automorphisms of a quantum
group, but the preceding remarks suggest something else.  Regarding $A(P)$ as a
quantum plane, one introduces a differential calculus ${\cal{D}}$ on $A(P)$.
The interesting object is the algebra of quantum automorphisms of this
differential calculus.

Here we shall obtain (Section 3), in the case that $A^*$ is the twisted, or
multiparameter, quantum group $U_{<q>,a}(gl_N)$, for generic $a$, the formula

$$U\hskip-1mm T = \prod_{{\scriptstyle i>j \atop \scriptstyle m<n} \atop
\scriptstyle k}
       e_a^{X_i^jP_i^j} e^{\tau_kH_k} e_{1/a}^{Y_m^nQ_m^n} \eqno(1.30)$$

\no The quantum group $A^*=U_{<q>,a}(gl_N)$ is generated by  $P_i^j$ (simple
positive roots), $Q_m^n$ (simple negative roots),  and $H_k$ (Cartan
generators); the dual    algebra \break $A(<q>,a)$  by $X_i^j$, $\tau_k$ and
$Y_m^n$;
$e_a^x$ is a deformed exponential.

The structure of  bialgebra on $A(A^*)$ determines the structure of  bialgebra
on
the dual, so that $U\hskip-1mm  T$ can be calculated along with the structure
of
$A^*(A)$.   The strategy followed in this paper was to begin from  the
structure
of the pseudo-group  $A$ = \A   (Woronowicz-Manin picture),  and use the
structural properties of the dual form to determine, first $U\hskip-1mm T$  and
then the structure of the quantum group  $A^* = U_{<q>,a}(gl_N)$ (Drinfeld
picture).  The  structural relations are, once more

$$T(x,l) \,\, T(x,l^\prime) = T(x,\Delta(l)), \eqno(1.31)$$

$$T(x,l) \,\, T(x^\prime,l) = T(\Delta(x),l) .\eqno(1.32)$$

\no Here $x$ and $x^\prime$ refers to two copies of $A$ and $l,l^\prime$ to two
copies of $A^*$; $\Delta(x)$ is the coproduct of $A$ and $\Delta(l)$ is the
coproduct of $A^*$.  The first relation determines the algebraic structure of
$A^*$ and the second one gives the coproduct (Section 4).

There are interesting homomorphisms into a quantum group from its dual, that
allow us to obtain the Universal R-matrix from the expression for the Universal
T-matrix.  If $\Phi$ is such a homomorphism, then

$$U\hskip-1mm R = (id \otimes \Phi ) U\hskip-1mm T; \eqno(1.33)$$

\no $\Phi$ operates on the generators of \A. (Section 6)

The similarity of (1.3) to (1.2) seems to suggest that the formula may be
useful
in connection with Fourier transforms on quantum groups.

\bb

\ce {\bf 2. Factorization}
\b

\no Let $R$ be the Yang-Baxter matrix (1.18), in the fundamental representation
of multiparameter $[12,13]$ (twisted $[11]$) quantum $gl(N)$, $V$ an
$N$-dimensional vector space, $P \in End(V \otimes V)$ defined by

$$P_{ji}^{kl} = R_{ij}^{kl} .\eqno(2.1)$$

\no Then $P$ satisfies the braid relation,

$$({\rm braid})_{123} \equiv P_{12}P_{23}P_{12} - P_{23}P_{12}P_{23} = 0,
\eqno(2.2)$$

\no and the Hecke condition (with {\it a} in the field )

$$(P-1)(P+a) = 0.\eqno(2.3)$$

Let $F$ be the algebra of formal power series finitely generated  by $(z_i^j),
\,\, i,j=1,\ldots,N$,
$(z_i^i)^{-1}, \,\, i=1,\ldots,N$, and the unit, with relations

$$z_i^i(z_i^i)^{-1} = (z_i^i)^{-1}z_i^i=1, \,\, i=1,\ldots,N.$$

\no The quantum algebra (pseudogroup [5]) $A =$ \A  is the quotient of $F$ by
the ideal  generated by the relations

$$[P,Z \otimes Z] = 0, \,\, Z \equiv {\rm matrix} \,\,(z_i^j).\eqno(2.4)$$

  There is a unique algebra homomorphism $\Delta: A \rightarrow A \otimes A$,
such that
$\Delta(Z) = Z \otimes T$; that is,

$$\Delta(z_i^j) = \sum_k z_i^k \otimes z_k^j, \eqno(2.5)$$

\no which gives $A$ a structure of bialgebra.

There is more than one sense in which \A is dual to the quantum group
$U_{<q>,a}(gl_N)$, a deformation of the enveloping algebra $U(gl_N)$ of
$gl(N)$.  Recall that (the differential of) the classical $r$-matrix defines  a
Lie structure on the dual space $gl(N)^*$, turning $\{gl(N), gl(N)^*\}$ into a
Lie bialgebra dual pair.  Here we shall relate \A to a deformation
$U_{<q>,a}(gl_N{}^*)$ of the enveloping algebra of $gl(N)^*$.  Then we show
that
there is a natural duality between the two deformed enveloping (bi-)algebras.

The first step is to justify the following factorization:

$$z_i^j = \sum_k X_i^k z_k Y_k^j, \eqno(2.6)$$

$$X_i^j = \cases{1, & $i=j$, \cr
                 0, & $i<j$, \cr} \quad
  Y_i^j = \cases{1, & $i=j$, \cr
                 0, & $i>j$, \cr} .\eqno(2.7)$$

\no Let $F^\prime$ be the algebra  generated by $X_i^j, \,\, z_k, \,\,
(z_k)^{-1}$ and $Y_i^j$ (with unit), with relations $z_k(z_k)^{-1} = (z_k)^{-1}
z_k=1$ and (2.7).  The formula (2.6) is invertible and allows the
identification
of $F^\prime$ with $F$.  The relations (2.4) make sense in
$F^\prime$ and the associated quotient is an alternative presentation of $A$.
We shall obtain the alternative expressions for the relations and show, in
particular, that the $X$'s commute with the $Y$'s.

\vskip.5cm
\no {\bf Proposition.} Let $A_-,A_+$ and $A_0$ be quotients of $A$ by the
ideals
generated by

$$I_- = \{z_i^j, \,\, i<j\}, \quad I_+ = \{z_i^j,\,\,i>j\}, \quad
        I_0 = \{z_i^j, \,\,i \not= j\}, \eqno(2.8)$$

\no then

$$A = A_- \mathrel{\mathop\otimes_{A_0}} A_+ .\eqno(2.9)$$

\vskip.5cm

\no {\bf Proof.}  The crucial ingredient is the fact that the sets $I_{\pm ,0}$
actually generate ideals of $A$.  This is an easy consequence of the relations,
but there is a deeper reason for it.  It is known that $A$ has representations
$\pi$ and $\pi^\prime$ given by

$$(\pi(z_i^k))_j^l = R_{ij}^{kl}, \quad (\pi^\prime(z_i^k))_j^l =
(R^{-1})_{ji}^{lk},\eqno(2.10)$$

\no reflecting the existence of two homomorphisms from $gl(N)^*$ to $gl(N)$.
The kernels of these two isomorphisms include the sets $I_-$ and $I_+$,
respectively.  Now let $(\tilde X_i^j)$ and $(\tilde Y_i^j)$ be the images  of
$(z_i^j)$ under the projections $A \rightarrow A_-$ and $A \rightarrow A_+$,
then
$\tilde X_i^j=0$ for $i<j$ and $\tilde Y_i^j=0$ for $i>j$, and  define

$$\tilde z_i^j \equiv \sum_k \tilde X_i^k \otimes \tilde Y_k^j.\eqno(2.11)$$

\no The relations among the $\tilde X$'s and among the $\tilde Y$'s are exactly
the same as the relations among the $z$'s and, because $\Delta$ is a
homomorphism, the same relations are also obeyed by the $\tilde z$'s.  The
composite mapping given by the two projections of $A$ on $A_+$ and on $A_-$
followed by (2.11) is one-one, so we can identify $\tilde z_i^j$ with $z_i^j$.
Now set

$$\tilde X_i^j=X_i^jx_j, \quad i \geq j, \quad \tilde Y_i^j = y_iY_i^j,
        \quad i \leq j, \quad  X_i^i=Y_i^i=1, \eqno(2.12)$$

\no and identify $X_i^j$ with $X_i^j \otimes 1, \,\, Y_i^j$ with $1 \otimes
Y_i^j$, so that

$$z_i^j = X_i^k z_k Y_k^j, \quad z_k \equiv x_k\otimes y_k .\eqno(2.13)$$

\no This is the desired formula (2.6). It is now evident that  the $X$'s
commute
with the $Y$'s, (2.13) is precisely what we mean (2.9), and the proposition is
proved.

As we said, the relations among the $\tilde X$'s and among the $\tilde Y$'s are
exactly the same as among the $z$'s, and one easily derives the following
relations among the generators of the decomposition (2.13).  First,

$$z_iz_j = z_jz_i,$$

$$z_kX_i^j = C_{kij}X_i^jx_k, \quad z_kY_i^j = C^\prime_{kij}Y_i^jz_k
.\eqno(2.14)$$

\no The coefficients $C_{kij}$ and $C^\prime_{kij}$ are given in the Appendix,
Eq. (A.3-3').  Next, define ``simple generators"

$$X_i \equiv X^i_{i+1}, \quad Y_i=Y_i^{i+1}, \quad i=1,\ldots,N,$$

\no then there are quommutation relations and Serre relations

$$[X_i,X_j]_{k_{ij}}=0, \quad [Y_j,Y_i]_{k_{ij}}=0, \quad |i-j|>1,
\eqno(2.15)$$

$$[X_{i+1},X_i]_{k_i}=(1-1/a)\,\,X_{i+1}^{i-1}, \eqno(2.16)$$

$$[Y_i,Y_{i+1}]_{k_i}=-(1-1/a)\,\,Y_{i+1}^{i-1}, \eqno(2.17)$$

$$[X_i,[X_{i+1},X_i]_{k_i}]_{r_i}=0=[X_{i+1},[X_{i+1},X_i]_{k_i}]_{s_i},
\eqno(2.18)$$

$$[[Y_i,Y_{i+1}]_{k_i},Y_i]_{r_i}=0=[[Y_i,Y_{i+1}]_{k_i},Y_{i+1}]_{s_i}.
\eqno(2.19)$$

\no Here $[A,B]_k \equiv AB-kBA$.  The coefficients $k_{ij}$ and
$k_i$  are in (A.4).

This alternative presentation of $A$, in terms of $X_i^j,Y_i^j,z_k$ and
$(z_k)^{-1}$, is very convenient for our purpose.  In particular, the
construction of the universal $T$-matrix for $A$ reduces to the same problem
for
$A_\pm$.  As an ultimate refinement, we introduce elements $\rho_k$ of
$A_-$, $\sigma_k$ of $A_+$ and $\tau_k$ of $A$ by setting

$$x_k=e^{\rho_k}, \quad y_k=e^{\sigma_k}, \quad z_k=e^{\tau_k}, \eqno(2.20)$$

\no and adopt $\rho_k, \sigma_k, \tau_k$ as generators instead of $x_k, y_k,
z_k$.  By abuse of notation we still use the same names, $A_\pm$ and $A$.  As
is
very well known, and obvious in the sequel, these algebras must be completed
with some infinite series, including the series (2.20).  Care must be taken to
extend far enough to get closure under the algebraic operations, without making
the algebras too large to be manageable.  For such questions we refer to the
papers [15] and [16].

In this new presentation, in which $(\tau_k)$ replace $(z_k^k)$ and
$(z_k^k)^{-1}$ as generators, $A$ becomes a deformation $U_{<q>,a}(gl_N{}^*)$
of
the enveloping algebra of the Lie algebra $gl(N)^*$ (with the Lie structure
determined by the classical $r$-matrix).  The Serre presentation of
$U_{<q>,a}(gl_N{}^*)$ is easily obtained from (2.14-19).

\bb

\ce {\bf 3. The Universal T-Matrix}
\b

\no The construction has already been explained elsewhere [10].   For $A_-$ we
take the basis

$$X^{[a][\alpha]} \equiv \prod_{\scriptstyle i>j \atop \scriptstyle k}
    (X_i^j)^{a_{ij}}(\rho_k)^{\alpha_k}, \quad e^{\rho_k} \equiv x_k
.\eqno(3.1)$$

\no To facilitate the manipulations that follow, it is crucial to adopt a good
ordering of the factors in the definition of the basis elements; the good rule
is that $X_i^j$  precedes $X_l^k$ if $j <k$ or $j=k, i < l$.

A rigorous definition of the dual is based on the natural basis that is
provided
by the functions defined by

$$P_{[a][\alpha]}(X^{[b][\beta]}) = \cases{1, &if
$[a][\alpha]=[b][\beta]$,\cr
                                          0, & \hbox{otherwise.}\cr}$$

\no For details concerning completion of the dual algebra in terms of entire
functions on $A$ we refer to the papers [15] by and [17].

If $\{P_{[a][\alpha]}\}$ is the dual basis, then the universal $T$-matrix for
$A_-$ is

$$T^- = \sum_{[a][\alpha]} X^{[a][\alpha]} P_{[a][\alpha]} . \eqno(3.2)$$

\no The important structural properties of this operator were discussed in the
Introduction; Eq. (1.32) for $T^-$ reads

$$\eqalign{\sum\bigl(\prod_{\scriptstyle i>j \atop \scriptstyle k}
          (X_i^j)^{a_{ij}}(\rho_k)^{\alpha_k}\bigr)
          &\bigl(\prod_{\scriptstyle i>j \atop \scriptstyle k}
          (X^{\prime j}_i)^{b_{ij}}(\rho^\prime_k)^{\beta_k}\bigr)
          P_{[a][\alpha]}P_{[b][\beta]}\cr
          &= \prod_{\scriptstyle i>j \atop \scriptstyle k}
          (\Delta X_i^j)^{c_{ij}}(\Delta\rho_k)^{\gamma_k}
          P_{[c][\gamma]} .\cr} \eqno(3.3)$$

\no We have abandoned the cumbersome notation with $\otimes$, writing $X$ for
$X \otimes 1$ and $X^\prime$ for $1 \otimes X$ from now on.  Comparing
coefficients of both sides one gets the relations sastisfied by the dual basis.
As an example, consider the coefficients of $(X_i^j)^{n-1} X^{\prime j}_i$ for
a
fixed pair $(i,j)$.  On the left one has (dots standing for zeros)
$P_{[\ldots,n-1,\ldots][0]} P_{[\ldots,1,\dots][0]}$; and on the other side the
only contributing term is

$$(\Delta X_i^j)^n P_{[\ldots,n,\dots][0]} . \quad (a_i^j=n, \quad
\hbox{all others zero}).$$

\no [This at first sight innocent statement is valid because of the particular
order chosen.] We need to know that the coefficient of $(X_i^j)^{n-1}X^{\prime
j}_i$ in
$(\Delta X_i^j)^n$ is

$$[n]_a \equiv (a^n-1)/(a-1),\eqno(3.4)$$

\no to get a simple recursion relation for $P_{[\ldots,n,\dots][0]}$, with the
solution

$$P_{[\ldots,n,\ldots][0]} = (P_i^j)^n/[n!]_a ,\eqno(3.5)$$

\no provided $a$ is not a root of 1.

The result is that the dual of $A_-$ is generated by

$$(P_i^j), \quad i>j, \quad {\rm and} \quad H_k=P_{[0][\ldots,1,\ldots]}, \quad
  \hbox{1 in k'th place}, $$

\no and that

$$T^- = \prod_{\scriptstyle i>j \atop \scriptstyle k}
        e_a^{X_i^jP_i^j} e^{\rho_kH_k}.$$

\no A similar calculation gives $T^+$, and the product $T^-T^+$, with the
identification $\rho_k+\sigma_k=\tau_k$ is the Universal $T$-matrix for
\A :

$$U\hskip-1mm T = \prod_{{\scriptstyle i>j \atop \scriptstyle m<n} \atop k}
      e_a^{X_i^jP_i^j} e^{\tau_kH_k} e_{1/a}^{Y_m^nQ_m^n} .\eqno(3.6)$$

\no The structure of $U\hskip-1mm T$ reflects that of (2.6) and (2.9).

\bbb

\ce {\bf 4. Relations and Coproduct of Twisted, Quantum $gl(N)$}
\b

\no Such relations are of course known [18], but we want to  show that they
drop
out of our formula for the universal $T$-matrix, and that the generators
$P_i^j$, $H_k$ and $Q_i^j$ are precisely the generators of a conventional
presentation.  Actually, we have found earlier derivations difficult to
understand, and we hope that the one given here may be an improvement.

First, in

$$T^-_{x,p} T^-_{x^\prime,p} = T^-_{\Delta x,p} \eqno(4.1)$$

\no we compare the coefficients of (properly ordered) elements of the basis of
$A$.  We first notice that terms linear in $X_i^j$ and $\rho^\prime_k$, for a
fixed triple $(i,j,k)$, occur with the same coefficients on both sides, as a
direct result of our construction.  But terms involving $\rho_kX^{\prime j}_i$
are in the wrong order; the coefficient on the left is $H_kP_i^j$, and on the
right side it is found by examination of the term

$$P_i^j \Delta(X_i^j) e^{H_k\Delta(\rho_k)} =
        P_i^j \Delta(X_i^j)(1+H_k\Delta(\rho_k) + \ldots).\eqno(4.2)$$

\no Applying (A.6-10) one sees that the relevant part of $\Delta(X_i^j)$ is
contained in

$$(x_i/x_j)X^{\prime j}_i = (1+\rho_i-\rho_j + \ldots)X^{\prime j}_i;$$

\no the coefficient we are looking for is thus
$P_i^j(H_k+\delta_{ik}-\delta_{jk})$, and

$$[H_k,P_i^j] = (\delta_{ki}-\delta_{kj})P_i^j.\eqno(4.3)$$

\no By these means one easily recovers the complete Serre presentation of the
positive Borel subalgebra of the quantum group $U_{<q>,a}(gl_N)$.  The
``simple"
generators are $P_i=P^i_{i+1}$, $i=1,\ldots,N$, and the remaining relations are

$$\eqalign{&[P_j,P_i]_{k_{ij}} = 0, \quad {\rm if} \quad |i-j|>1,\cr
           &[[P_i,P_{i+1}]_{k_i},P_i]_{r_i} =
          [[P_i,P_{i+1}]_{k_i},P_{i+1}]_{s_i}=0 .\cr} \eqno(4.4)$$

Next, we get the commutator between the simple generators of
$U_{<q>,a}(gl_N)$,

$$P_i \equiv P^i_{i+1} \quad \hbox{and} \quad Q_i \equiv Q_i^{i+1},
\eqno(4.5)$$

\no by comparing coefficients of $X_{i+1}^{\prime i}Y^{i+1}_i$ on both sides of
$T_{x,p} T_{x^\prime,p} = T_{\Delta x,p} $. For this we need to know parts of
the expression for $\Delta(\rho_k)$.  This is found in the Appendix; the result
is that, for $i = 1,2,...,N-1,$

$$[P_i,Q_j] = \delta_i^j{a \over 1-a} (q^{i,i+1})^{1-H_{i+1}-H_i}(a^{-H_{i+1}}
               - a^{-H_i})C_i.\eqno(4.6)$$

\no The relations

$$\eqalign{&[H_k,Q_i^j] = (\delta_{ki}-\delta_{kj}) Q_i^j  \cr
           &[Q_i,Q_j]_{k_{ij}} = 0, \quad \hbox{if} \quad |i-j| = 1, \cr
           &[Q_i,[Q_{i+1,},Q_i]_{k_i}]_{s_i} = [Q_{i+1,}
            [Q_{i+1,}Q_i]_{k_i}]_{s_i} = 0, \cr} \eqno(4.7)$$

\no complete the structure.  The dependence on $q^{i,i+1}$ can be removed by a
slight redefinition of the simple generators to make the structure reduce to
that of standard quantum $gl(N)$, as found in [18], but that would mess up the
expression (3.6).

The coproduct is obtained from the other structural formula, Eq. (1.4),

$$  \prod_{{\scriptstyle i>j \atop
           \scriptstyle m<n} \atop \scriptstyle k}
           e_a^{X_i^jP_i^j} e^{\rho_k H_k} e^{Y_m^nQ_m^n}
           \prod_{{\scriptstyle i>j \atop
           \scriptstyle m<n} \atop \scriptstyle k}
           e_a^{X_i^jP^{\prime j}_i} e^{\rho_k H^\prime_k}
           e^{Y_m^n Q^{\prime n}_m}
           = \prod_{{\scriptstyle i>j \atop
           \scriptstyle m<n} \atop \scriptstyle k}
           e_a^{X_i^j\Delta P_i^j} e^{\rho_k \Delta H_k}
           e^{Y_m^n \Delta Q_m^n} ,  \eqno(4.8)$$

\no it is the homomorphism generated by

$$\Delta(H_k) = H_k \otimes 1 + 1 \otimes H_k,$$
\vskip-5mm
$$\Delta(P_i) = P_i \otimes 1 + A_i \otimes P_i, \quad
  \Delta(Q_i) = Q_i \otimes B_i + 1 \otimes Q_i .$$
\vskip-5mm
$$A_i = C_i(q^{i+1,i})^{H_i+H_{i+1}} a^{-H_i}, \quad
  B_i = C_i(q^{i+1,i})^{H^\prime_i+H^\prime_{i+1}} a^{-H^\prime_{i+1}},$$
\vskip-5mm
$$C_i = \prod_{k\not= i,i+1} (q^{i+1,k}q^{ki})^{H_k}.\eqno(4.9)$$

\bb

\ce {\bf 5. Applications of the Universal T-Matrix}
\b

\no We noted the existence of two representations of the quantum algebra \A in
$gl(N)$.  If, in the formula (3.6), one takes the generators $P, H$ and $Q$ in
the fundamental representation, then one recovers the original
$N$-dimensional matrix $Z=(z_i^j)$, and if one then takes $z_i^j$ in the
representation (2.10), then one gets the original $R$-matrix in the fundamental
representation.  Our results concerning the structure of the dual algebra shows
that the representations (2.10) [see Appendix, Eq.(A.9)] lift to the
structure;
there are a homomorphism $\Phi, \Phi^\prime: A \rightarrow U_{<q>,a}(gl_N)$
given
by

$$\Phi(X_j^i)=(1/a-1)q^{ji} Q_i^j, \quad i<j, \quad \Phi (Y_i^j) = 0, $$

$$\Phi(z_k) = \bigl(\prod_i (q^{ki})^{H_i}\bigr) a^{H_{k+1}+\ldots+H_N},
\eqno(5.1)$$

$$\Phi^\prime(Y_i^j)=(a-1)q^{ji} P_j^i, \quad i<j, \quad \Phi'(X_i^j) = 0, $$

$$\Phi^\prime(z_k)=\bigl(\prod_i
(q^{ki})^{H_i}\bigr)a^{-H_1-\ldots-H_{k-1}}.\eqno(5.2)$$

\no One can view the universal T-matrix as an element of
$U_{<q>,a}(gl_N{}^*) \otimes U_{<q>,a}(gl_N)$.  Then we have

\vskip.5cm

\no {\bf Theorem.}  The universal R-matrix for $U_{<q>,a}(gl_N)$ is found by
applying the mapping $id \otimes \Phi $ to the universal $T$-matrix (3.6),

$$U\hskip-1mm R=(id \otimes \Phi ) U\hskip-1mm T =
\prod_{\scriptstyle i>j \atop \scriptstyle k}
     e_a^{\Phi(X_i^j)P_i^j} [\Phi(z_k)]^{H_k}. \eqno(5.3)$$

%\no {\bf Proof.}  Universal R-matrices for twisted quantum $gl(N)$ are known
%[19]; by inspection one verifies that mappings with the structure of $\Phi$
%%and
%$\Phi^\prime$ exist; that they are homomorphisms is a consequence of the
%Yang-Baxter equation.  To demonstrate that they take the form (5.1) and (5.2)
%it is evidently enough to examine their restriction to the simple roots.  It
%should be very instructive to have a more direct proof.

\bb

\ce {\bf 6. Further Generalizations}
\b

\no 6.1. STANDARD QUANTUM GROUPS

These are the 1-parameter deformations given by Drinfeld in [1]. Much of their
structure is summed up in the quantum double construction; it is basically the
same as for twisted $gl(N).$ There are representations defined by (2.10), with
kernels $I_{\pm}$ and associated quotients $A_{\pm}$, and the structure formula
(2.9),

$$A = A_- \otimes_{A_0} A_+  $$

\no holds generally, with $A_0$ the Abelian quotient by $I_+ \cup I_-$. As for
twisted $gl(N)$, these representations lift to isomorphisms $\Phi$  and $\Phi'
$
from $A_{\pm }$ to subalgebras $U_{\pm }$ of $U$. The  Universal R-matrix is
known [19] (hence, so  are $\Phi$ and $\Phi'$),  actually in two forms, $R_+
\in
U_- \otimes U_+$  and $R_- \in U_+ \otimes U_-$. The Universal T-matrices for
$A_{\pm}$  are

$$T^- = (\Phi^{-1} \otimes 1)R_-, \quad T^+ = (\Phi'^{-1} \otimes 1)R_- $$

\no and for $A$, it is $$U\hskip-1mm T = T^-T^+.$$

\no 6.2. NONSTANDARD QUANTUM GROUPS

These include the constructions on $A_n^{(2)}$ and $D_n^{(2)}$ of Jimbo [2], as
well as ``Esoteric Quantum $gl(N)$" and the other deformations of twisted
quantum
$gl(N)$. It seems that the general formulas that apply to the standard quantum
groups have direct generalizations to these cases as well. However, the
structure of the quantum double is very different and $A_0$ is no longer
abelian. See [10], second paper.

\vskip1cm

\no 6.3. ROOTS OF UNITY

The most interesting aspects of quantum groups are those that appear at special
values of the parameters; up to this point, in order to postpone the discussion
of these phenomena, we have assumed $(< \hskip-1mm q
\hskip-1mm>,a)$ in generic position.

Recall that the exponential form $U\hskip-1mm T = e^{xp} $  appeared in our
first example, in subsection (1.3), as the series $\sum F_kx^k$ with $F_k$
subject to the recursion relation (1.10),

$$F_kF_1 = (k+1)F_{k+1} \quad \Rightarrow \quad F_k = p^k/k!. \eqno(6.1)$$

\no In the evaluation of $U\hskip-1mm T$ for  the quantum groups one encounters
modified exponentials,

$$e_a^{xp} = \sum(xp)^n/[n!]_a, \eqno(6.2)$$

\no arising from the solution of the recursion relation

$$F_kF_1 = [k+1]_aF_{k+1} \quad \Rightarrow \quad F_k = p^k/[k!]_a.
\eqno(6.3)$$

\no See Eq.(3.5). This is the only place that the restriction to generic
parameters is relevant. There are no restrictions on the parameters $q^{ij}$.
The limitation to parameters in generic position applies only to the Hecke
parameter $a$, and the (till now) excluded values are those for which there is
an integer $K$ such that

$$a^K = 1, \quad a^n \neq 1, \quad n = 1,2,...,K-1.$$

\no We now suppose that this relation holds; thus $a$ is a primitive root of
unity.

Then $[K]_a = 0,$ and the general solution of the recursion relation (6.3) is

$$F_{mK+n} = (p'^m/m!)(p^n/[n!]_a), \quad p^K = 0,$$

\vskip-5mm

$$m = 0,1,... ; \quad n = 0,1,..., K-1.$$

\no This involves a new generator $p'$, and the new relation $p^K = 0$. The
expressions obtained for the Universal T-matrices are thus modified, when
$a$ is a root of unity, by the replacement of all twisted exponentials,
according to the rule

$$e_a^{xp} \rightarrow \sum_m (x^Kp')^m/m!\sum_{n=0}^{K-1}(xp)^n/[n!]_a,$$

\no with $p^K = 0$.

What is significant is the appearance of a new independent generator $p'$ that
replaces $p^K$. The structure of the dual algebra is drastically changed; it is
no longer the Drinfeld quantum group $U_{<q>,a}(\cal G)$.  Of course, this
latter still exists at roots of unity, but its dual is not the pseudogroup \A.

The structure of the dual of \A at $a^K = 1$ can be found by the same methods
as
in the generic case, but much more quickly by the prescription

$$p' := \lim_{a^K \rightarrow 1} p^K/[K]_a.$$

\no Thus one obtains for the dual of \A the additional relations

$$ [P_i, P'_j] = 0, \quad [P'_i, P'_j] = 0,$$
\vskip-5mm
$$[H_k, P'_i{}^j] = K(\delta_{ki} - \delta_{kj})P'_i{}^j, $$
\vskip-5mm
$$[P_i,Q'_i] =(a-1)(q^{12})^{1-H_1-H_2}\left(Q_1^{K-1}a^{-H_2} -
a^{-H_1}Q_1^{K-1}\right).$$

\no The last relation is noteworthy; it shows that the  generators $P'_i$ {\bf
do not} generate an ideal; consequently, the  irreducible modules  {\bf do not}
become nondecomposable at roots of unity.

The dual of $U_{<q>,a}(gl_N)$ can be obtained in the same way. The result is
that $(X_i^j)^K$ and $(Y_i^j)^K$ vanish and new independent generators
$X'_i{}^j$ and $Y'_i{}^j$ appear. In this case $X'_i{}^j$ and $Y'_i{}^j$ 
do generate ideals.

In the simplest  case of $U_{q,q'}(gl_2)$ the dual  pseudogroup is  generated
by
$X,X',\break \rho_1,\rho_2,Y$ and $Y'$, with relations

$$X^K = Y^K = 0, \quad q = e^h, \,\,q' = e^{h'},$$
\vskip-5mm
$$[\rho_1,X] = hX,\,\, [\rho_1,X'] = KhX',$$
\vskip-5mm
$$ [\rho_2,X] = h'X ,\,\, [\rho_2,X'] = Kh'X', etc. $$

\no The associated Woronowicz matrix; that is, $U\hskip-1mm T$  evaluated in
the
2-dimensional representation of $gl(2)$, is

$$\pmatrix{1 & 0 \cr X & 1 \cr} \pmatrix{ e^{\rho_1} & 0 \cr 0 & e^{\rho_2}
\cr}     \pmatrix{1 & Y \cr 0 & 1 \cr},$$
\vskip3mm
\no in which $X'$ and $Y'$ do not appear. To obtain the full dual  one has to
take $U\hskip-1mm T$ in a representation of quantum $gl(N)$ in which $P^K \neq
0$.

Duality at roots of unity has already been described by Fr\"olich  and Kerler
[20], for the case of standard quantum $sl(2)$.

Classical q-functions are always investigated in the domain $\mid q \mid  < 1$.
Indeed, many of these functions, including the twisted exponential, cease to
exist as $q$ tends to a root of unity. We have seen that the twisted
exponential, within the algebraic contest, does have a natural generalization
to
roots of unity: but that this generalization is highly nontrivial and even
mysterious in a purely analytical  context. It is proposed to take up the study
of $q$-functions on a broader basis, with a view to finding natural extensions
to roots of unity. Perhaps even the Rogers-Ramanujan identies can be
generalized
to $q^K = 1$.

\bb

\ce {\bf 7. Quantum $gl(N)$ and Quantum $sl(N)$.}

\vskip5mm

\no The (Drinfeld) quantum group $U_{<q>,a}(sl_N)$ is the subquotient of
$U_{<q>,a}(gl_N)$ defined by the ideal generated by the element

$$\sum_{i+1}^N H_i \equiv \cal H.$$

\no The universal T- and R-matrices of quantum $gl(N)$ reduce to those of
quantum $sl(N)$ under the projection that anulls $\cal H$. We shall calculate
this reduced R-matrix, in the fundamental representation  of $sl(N)$.

Restricting the Universal R-matrix  (5.3) to any faithful N-dimensional
representation gives the formula

$$ R = \left[ 1 + \sum_{i>j}P_i^j\Phi (X_i^j) \right]
\sum_{i,j}(\tilde q^{ij})^ {H_i \oplus H_j}, \eqno(7.1) $$

\no with $\tilde q^{ij} = a q^{ij}$ for $i<j $ and $\tilde  q^{ij} = q^{ij}$
for
$i \ge j$. With ${\cal H} = 0$, in the $N$-dimensional representation of
$sl(N)$
we have $H_i = M_i^i - 1/N$. Substituting this into (7.1) we obtain

$$R = (1-a)\sum_{i<j} (\kappa_i/\kappa_j M_j^i \otimes M_i^j  + \sum_{i,j} \hat
q^{ij} M_i^i \otimes M_j^j, \eqno(7.2)$$

\no with

$$\hat q^{ij} = (\kappa_i/\kappa_j)q^{ij}, \quad
\kappa_i = \left( a^i \prod_k q^{ki} \right)^{1/N}. \eqno(7.3) $$

\no This R-matrix for $sl(N)$ differs from that of $gl(N)$ in two particulars.
(1) $q^{ij}$ is replaced by $\hat q^{ij}$. These new parameters are not
independent:

$$\prod_i \hat q^{ij} a^j = a^{(N+1)/2}. \eqno(7.4)$$

\no If the original parameters satisfy

$$\prod_i q^{ij} a^j = a^{(N+1)/2}, \eqno(7.5)$$

\no then $\hat q^{ij} = q^{ij}$. Therefore, reducing from $gl(N)$ to $sl(N)$ as
we have done, from arbitrary initial parameters, gives the same result as
restricting the parameters to satisfy (7.5).  (2) The factor
$\kappa_i/\kappa_j$
can be removed by the isomorphism $M_i^j \rightarrow (\kappa_i/\kappa_j)M_i^j$.

The restriction (7.5) was found by Schirrmacher [12]  from the requirement that
the quantum determinant be unity.

The deformations of \A have recently been calculated [15]. This algebra is
rigid
for essential, first order deformations except for very special values of the
parameters. There are several series of deformations of twisted quantum
$gl(N)$, here we illustrate the simplest one, in the case of $gl(3)$.

The parameters are $q^{12}, q^{23}, q^{13}$ and $a$. The deformation that
consists of adding the following piece

$$\delta R = \epsilon \left( q^{13} M_1^2 \otimes M_3^2  -  M_3^2 \otimes M_1^2
\right) $$

\no to (1.18), is exact and essential. It exists if and only if the parameters
satisfy

$$ q^{12} = q^{23} \quad {\rm and} \quad q^{13} = (q^{12})^2.$$

\no The projection to $sl(3)$ fixes the value of $a$:

$$ q^{12} = q^{23} :=q, \quad q^{13} = q^2. \quad a = q^{-3}.$$

That is why, for $sl(3)$, unlike $gl(3)$, the term  ``roots of unity" applies
to
the values of the $q'$s.

This esoteric form of quantum $sl(3)$ is not included in the list produced by
Jimbo [2], but the classical limit is in the classification of  Belavin and
Drinfeld [21].

\bbb

\ce {\bf Appendix}
\b

\no The $R$-matrix for twisted, quantum $gl(N)$, in the fundamental
representation, was given in (1.18),

\vskip-5mm

$$ R = \sum_i M_i^i \otimes M_i^i
             + \sum_{i<j}
           \bigl(q^{ji}M_j^j \otimes M_i^i + aq^{ij}M_i^i
           \otimes M_j^j + (1-a)M_j^i \otimes M_i^j \bigr), $$

\no   The inverse matrix is given by the same formula, with the parameters
$q^{kl}$ and $a$ replaced by their inverses.

The relations (2.4) of \A, written out in full detail, are

\vskip-5mm

$$\eqalignno{z_i^az_i^b &=q^{ab}z_i^b z_i^a, \cr
             z_i^az_j^a &=(aq^{ij})^{-1}z_j^az_i^a, \,\, i<j, \cr
     z_i^az_j^b &=(aq^{ab}/q^{ij})z_j^bz_i^a, \,\, i>j, \,\, a<b \cr
     q^{ij}z_i^az_j^b &-q^{ab}z_j^bz_i^a=(a-1)z_j^az_i^b, \,\, i>j, \,\,
a>b.&(A.1)\cr}$$

\no The quotient algebras $A_+,A_-$ satisfy relations that one gets from these
by annulling the generators in $I_+,I_-$,

$$A_-: z_i^j \rightarrow \cases{\tilde X_i^j, & $j \leq i$ \cr
                                 0, & $i<j$ \cr}; \quad
  A_+: z_i^j \rightarrow \cases{\tilde Y_i^j, & $i \leq j$ \cr
                                 0, & $j<i$ \cr} \eqno(A.2)$$

\no Finally, one gets from (2.12) the following relations for $A_-$,

\vskip-5mm
$$x_ix_j = x_jx_i, $$
\vskip-5mm
$$x_kX_i^j=\cases{q^{ik}q^{kj}X_i^jx_k,
\hfill& $k<j<i$ \hbox{or} $j<i\leq k$,\cr
      (1/a)q^{ik}q^{kj}X_i^jx_k, \hfill& $j \leq k<i$\cr} \eqno(A.3)$$
\vskip2mm
\no These become relations of $A$ after substitution of $z_i$ for $x_i$.
Likewise, the following relations for $A^+$ become relations of $A$ if $y_i$ is
replaced by $z_i$,
\vskip-5mm
$$y_iy_j=y_jy_i,$$
\vskip-5mm
$$y_kY_i^j=\cases{q^{ik}q^{kj}Y_i^jy_k,\hfill& $k \leq i<j$ \hbox{or}
$i<j<k$,\cr
          aq^{ik}q^{kj}Y_i^jy_k, \hfill& $i<k \leq j$. \cr} \eqno(A.3')$$
\vskip2mm
\no The coefficients in (2.15-19) are

$$k_{ij}=q^{i+1,j}q^{ji}/q^{i+1,j+1}q^{j+1,i},$$
\vskip-5mm
$$k_i=q^{i+1,i-1}/q^{i+1,i}q^{i,i-1}.\eqno(A.4)$$

We need some formulas for the coproduct of $A_-$.  From (2.12) and

$$\Delta(\tilde X_i^j) = \sum_k \tilde X_i^k \otimes \tilde X_k^j \eqno(A.5)$$

\no one gets
\vskip-4mm
$$\Delta(x_i) = x_i \otimes x_i, \eqno(A.6)$$
\vskip-8mm
$$\Delta(X_i^j)=\sum_k X_i^k(x_k/x_j) \otimes X_k^j=X_i^j \otimes
         1+(x_i/x_j) \otimes X_i^j + \ldots .\eqno(A.7)$$

\no Only the first two terms are relevant for our calculation of $T^-$;
neglecting the rest we have

$$\Delta(X) = A+B, \quad A=X_i^j \otimes 1, \quad B=(x_i/x_j) \otimes
X_i^j.\eqno(A.8)$$

\no One has $BA=aBA$ from (A.3) and thus
$(A+B)^n=A^n+A^{n-1}B(1+a+\ldots+a^{n-1})$ and finally the result

$$\Delta(X_i^j)^n=(X_i^j)^n \otimes 1+[n]_a(X_i^j)^{n-1}(x_i/x_j)
  \otimes X_i^j + \ldots ,$$

\no that was used in (3.4).

For the evaluation of $[P_i,Q_i]$ we need to know some terms in $\Delta(z_i)$,

$$\Delta(z_k) = (z_kz_k^\prime)^{H_k}\biggl(1+(q^{k,k+1}/a)Y_kX^\prime_k
                 -q^{k-1,k}X^\prime_{k-1}Y_{k-1}\ldots\biggr)$$
\vskip-8mm
$$\bigl(\Delta(z_k)\bigr)^{H_k}=(zz^\prime)^{H_k}
   \bigl(1+[H_k]_a(q^{k,k+1}/a)Y_kX^\prime_k -[H_k]_{1/a}
   q^{k-1,k}Y_{k-1}X^\prime_{k,k-1}+\ldots \bigr).$$

\no The relation (4.6) now follows easily in the same way.

The representations $\pi$ and $\pi^\prime$ are given by (2.10) and (1.18), and
more explicitly by vskip-3mm
$$\pi(z_j^i)=\cases{(1-a)M_i^j , \hfill& $i<j$,\cr
                    0\quad , \hfill& $i>j$,\cr} \quad
\pi(z_k^k)_i^j = \delta_i^j q^{ki} a^{(i>k)},$$

$$\pi^\prime(z_i^j)=\cases{(1-1/a)M_j^i, \hfill& $i<j$,\cr
                            0, \hfill& $i>j$,\cr} \quad
\pi^\prime(z_k^k)_i^j = \delta_i^j q^{ki} a^{-(i<k)},\eqno(A.9)$$

\no where $(i<k)=1$ if $i<k$ and zero otherwise.

\bbb

\ce {\bf Acknowledgments}

This work was supported in part by the National Science Foundation.

\bbb

\ce {\bf References }

\no  1. V.G. Drinfeld, Int.Congr.Math., Berkeley 1986 798-820.

\no  2. M.Jimbo, Commun.Math.Phys. {\bf 102} (1986) 537-547.

\no  3. R.J. Baxter, Ann.Phys. {\bf 70} (1971) 193-228.

\no  4. E.K. Sklyanin and L.D. Faddeev, Sov.Phys.Dokl. {\bf 23} (1978) 902-904.

\no  5. L. Woronowicz, Publ.Res.Ins.Math.Sci., Kyoto Univ.,  {\bf 23} (1987)
117-181.

\no  6. Y. Manin, {\it Topics in noncommutive geometry},  Princeton Univ.
Press,
1991.

\no  7. Y.N. Reshetikhin, L. Takhtajan and L.D. Faddeev, Leningrad Math.J.
        {\bf 1} (1998)
\line{\hfill 193-225;  A. Sudbery, J.Phys.A {\bf 20 } (1990), L697.}

\no 8. P. Bonneau, M. Flato and G. Pinczon, Lett.Math.Phys. {\ bf 25} (1992)
 75-84.

\no  \line{9. S.P. Vokos, B. Zumino and J. Wess, Z.Phys.C  {\bf 48} (1990)
65-74; \hfil}
\line{\hfil C. Fronsdal, Lett.Math.Phys. {\bf 24} (1992) 73-78.}

\no 10. C. Fronsdal and A. Galindo, The dual of a
   Quantum Group,  Lett.Math.Phys.
 {\bf 27} (1993)  59-71;  The Universal T-Matrix,  Proc. Joint Summer Res.Conf.
\break  (AMS-IMS-SIAM) on Conformal Field Theory,  Topological Field Theory and
\line{\hfill Quantum Groups.  [Preprint  UCLA/93/TEP/ 2]}

\no 11. N.Y. Reshetikhin, Lett.Math.Phys. {\bf 20} (1990) 331-336.

\no 12. A. Schirrmacher, Z.Phys.C {\bf 50} (1991) 321.

\no 13. A. Sudbery, J.Phys.A {\bf 23} (1990) L697.

\no 14. J.C. Baez, Lett.Math.Phys. {\bf 23} (1991) 133.

\no 15. C.Fronsdal and A.Galindo, Deformations of Multiparameter quantum
$gl(N).$
\line{\hfill Preprint UCLA/92/TEP/52.}

\no 16. P. Truini and V.S. Varadarajan, Lett.Math.Phys. {\bf 21} (1991) 287.

\no 17. M.V. Pimsner, A Class of Markov Traces, preprint.

\no 18. J.F. Cornwell, J.Math.Phys. {\bf 33} (1992) 3963-3977).

\no 19. Kirillov and N.Y. Reshetikhin, Adv.Series in Math.Phys. {\bf 7} (1988)
285-339.

\no 20. J. Fr\"olich and T.Kerler, {\it Quantum Groups, Quantum Categories and
Quantum }
\line{\hfill {\it  Field Theory}, to be published.}

\no 21 A.A. Belavin and V.G. Drinfeld, Sov.Sci.Rev.Math. {\bf 4} (1984) 93-165.

%\vskip.3cm
%\no C. Fronsdal

%\no{\it Department of Physics}

%\no {\it University of California, Los Angeles, CA 90024-1547}

\bye